\documentclass[12pt]{amsart}
\usepackage[english]{babel}
\usepackage{amsmath}
\usepackage{amsthm}
\usepackage{amssymb}
\usepackage{mathrsfs}
\usepackage{enumerate}
\usepackage[notcite, final, notref]{showkeys}
\usepackage{amsfonts}
\usepackage{dsfont}
\usepackage{tikz}
\usepackage[siunitx]{circuitikz}
\usetikzlibrary{arrows,calc,chains,shapes,dsp}
\def\C{\mathbb C}
\def\R{{\mathbb R}}
\newtheorem{Pa}{Paper}[section]
\newtheorem{Tm}[Pa]{{\bf Theorem}}
\newtheorem{La}[Pa]{{\bf Lemma}}

\newtheorem{Cy}[Pa]{{\bf Corollary}}
\newtheorem{Rk}[Pa]{{\bf Remark}}
\newtheorem{Pn}[Pa]{{\bf Proposition}}

\newtheorem{Ex}[Pa]{{\bf Example}}
\newtheorem{Dn}[Pa]{{\bf Definition}}
\usepackage{xcolor}

\title[Passive Systems \& Matrix-convex invertible sets 
]
{Passive Linear Discrete-time Systems:\\[0.2cm]
Characterization through Structure}
\author[I. Lewkowicz]{Izchak Lewkowicz}
\address{School of Electrical and Computer Engineering 
Ben-Gurion University of the Negev\\ P.O.B. 653\\ Beer-Sheva, 84105\\
Israel}
\email{izchak@bgu.ac.il}

\begin{document}

\begin{abstract}
We here show that the family of finite-dimensional, discrete-time,
passive, linear time-invariant systems can be characterized through
the structure of a matrix-convex set, which is maximal in the sense of
being closed under products of its elements Moreover, this
observation unifies three setups:
(i) difference inclusions, \mbox{(ii) matrix-valued} rational functions,
(iii) realization arrays associated with rational functions.
\vskip 0.2cm

\noindent
It turns out that in the continuous-time case the corresponding 
structure is of a maximal matrix-convex cone closed under inversion.
\end{abstract}
\maketitle


\date{today}
\tableofcontents

\bibliographystyle{plain}
\section{Introduction}
\setcounter{equation}{0}

\noindent
In the study of dynamical systems, passivity is a fundamental property.
Thus, it has been extensively addressed in various frameworks.
A fundamental contribution was made by J.C. Willems, see e.g.
\cite{Willems1st1972}, \cite{Willems2nd1972} and \cite{Willems1976}.
Here we confine the discussion to discrete-time systems. More precisely,
we focus on the finite-dimensional, linear time-invariant case. These
passive systems are modelled by {\em Discrete-time Bounded} real
rational functions, denoted by $~\mathcal{DB}$, namely (here, for
simplicity, in a scalar framework) rational functions which map
\mbox{$\{z\in\C : |z|>1\}$}, the exterior of the closed unit disk, to
\mbox{$\{z\in\C : 1\geq |z|\}$}, the closed unit disk. (In the sequel,
the discussion is of matrix-valued rational functions).
\vskip 0.2cm

\noindent
In reading relevant literature, some extra care is quite helpful:\\
Here we follow the engineering motivation where $F(z)$ can be interpreted
as the $Z$-transform, i.e. \mbox{$F(z)=C(zI_n-A)^{-1}B+D$}, of the
shift-invariant difference state-equation $x(k+1)=Ax(k)+Bu(k)$,
$y(k)=Cx(k)+Du(k)$,
where in some sense, the input $u$ dominates the output $y$, i.e.
$\|u\|\geq\|y\|$ for all input $u$. For more general formulation
see e.g. \cite{GoeVanGesSuyVanDoorenDeMoor2003}.
\vskip 0.2cm

\noindent
In contrast, apparently motivated by symmetry, mathematical-analysis
circles prefer to study {\em Schur} functions which analytically map
the open unit disk \mbox{$\{z\in\C~:~1>|z|~\}$} to its closure,
\mbox{$\{z\in\C~:~1\geq |z|~\}$}, see e.g. \cite{deBraRovn1966} and
\cite{Schur}. This family was for example addressed in
\cite{AlpDijkRovndeSnoo} and \cite{BallStaffans2006}.
\vskip 0.2cm

\noindent
On the top of this difficulty, there is a confusion in names of
various associated families of functions. See e.g.
\cite[Section 2]{AlpayLew2011}, \cite{Lewk2020b} and
\cite{LiuXio2018}. We shall not further pursue this point.
\vskip 0.2cm

\noindent
There have been characterizations of finite-dimensional discrete-time
passive systems. For a modest sample of the literature on the subject,
see e.g.  \cite{Najs2013}, \cite{PremJury1994}, \cite{Vaidya1985} and
\cite[Section 4]{XiaoHill1999}. For an account
of the infinite-dimensional case, see e.g. \cite{ArliHassdeSnoo2007},
\cite{ArovStaff2005}, \cite{BallGroeterH2018}, \cite{BallStaffans2006}
and \cite{Staffans2006}.
\vskip 0.2cm

\noindent
Here we adopt a more abstract point of view and focus on the following
question:
\begin{quote}
How can one characterize the family of finite-dimensional,
discrete-time, passive, linear time-invariant systems 
through the structure of the whole set?
\end{quote}
The answer is that this family forms a {\em matrix-convex} set which
is maximal with respect to being closed under products among its
elements. Moreover, this observation unifies three setups:
\vskip 0.2cm

\begin{itemize}

\item[(i)~~]{}Difference inclusions,

\item[(ii)~]{}Discrete-time Bounded real rational functions,

\item[(iii)]{}Families of realization arrays of $\mathcal{DB}$, 
Discrete-time Bounded real rational functions.
\end{itemize}
\vskip 0.2cm

\noindent
This is a follow-up of the study in \cite{Lewk2020a} of the continuous-time
case. Combining the message of both works can be summarized as follows,
\vskip 0.2cm

\begin{center}
\begin{tabular}{|c|c|}
\hline
\multicolumn{2}{|c|}
{\rm Passive linear time-invariant systems and matrix-convexity}
\\
\hline
{\rm discrete-time}&{\rm continuous-time}\\
\hline
{\rm a maximal set closed under}&
{\rm a cone closed under inversion and}\\
{\rm products of its elements}&
{\rm maximal non-singular/analytic}\\
\hline
\end{tabular}
\end{center}

\noindent
This work is organized as follows. In Section \ref{Sec:DiffIncl} we
lay the foundation to the sequel and recall in sets of matrices all
satisfying a Stein inclusion with the same factor. Then, in Section
\ref{Sec:MCsets} we restrict the discussion to matrix-convex sets
of matrices which are closed under products among its elements and
maximal in this sense,
see Proposition \ref{Pn:MCIC}. As a sample motivation we recall in
the problem of stability of difference inclusions.
\vskip 0.2cm

\noindent
Subsequently, in Section \ref{Sec:MCsetsOfSystems} we exploit
Proposition \ref{Pn:MCIC} to characterize, in Proposition
\ref{Ob:DTBR}, Discrete-time Bounded, $\mathcal{DB}$, real rational
functions. Then the same structure is used to describe, in 
Corollary \ref{CyConvexRealizationS}, families of realization
arrays, of $\mathcal{DB}$ functions.
\vskip 0.2cm

\noindent
The conclusion that passive discrete-time systems (even
of various dimensions) are inter-related, is illustrated in
Examples \ref{Ex:MatrixConvexRational} and
\ref{Ex:MultiplicativeRealizationArrays}, for rational functions
and for realization arrays, respectively.

\section{Sets of Matrices with Common Stein Factor}
\label{Sec:DiffIncl}
\setcounter{equation}{0}

We start with notations. Let $\overline{\mathbf H}_n$ (${\mathbf H}_n$)
be the set of $n\times n$ Hermitian (non-singular) matrices and by
($\overline{\mathbf P}_n$) ${\mathbf P}_n$ denote the subsets of
$n\times n$ positive (semi)-definite matrices.
\vskip 0.2cm

\noindent
Now, for a prescribed $H\in\mathbf{H}_n$, consider the set
of all $n\times n$ matrices $A$ sharing the same Stein factor,
\begin{equation}\label{eq:DefSteinH}
\begin{matrix}\mathbf{S}_H&=\left\{ A\in\C^{n\times n}:~
H-A^*HA\in\mathbf{P}_n\right\}
\\~\\
\overline{\mathbf{S}}_H&=\left\{ A\in\C^{n\times n}:~
H-A^*HA\in\overline{\mathbf P}_n\right\}\end{matrix}
\quad\quad\quad H\in\mathbf{H}_n~.
\end{equation}
The set $\overline{\mathbf{S}}_H$ is the closure of the open set
$\mathbf{S}_H$ in the sense that $\overline{\mathbf P}_n$ is the
closure in $\overline{\mathbf H}_n$ of the open set $\mathbf{P}_n$.
\vskip 0.2cm

\noindent
One can refine the above definition by adding a parameter $\alpha>0$ to
obtain,
$\overline{\mathbf{S}}_H$),
\begin{equation}\label{eq:DefRefinedSteinH}
\begin{matrix}
{\scriptstyle\frac{1}{\alpha}}\mathbf{S}_H=\left\{ A\in\C^{n\times n}:~
H-{\scriptstyle\frac{1}{{\alpha}^2}}A^*HA\in\mathbf{P}_n\right\}
\\~\\
{\scriptstyle\frac{1}{\alpha}}\overline{\mathbf S}_H=\left\{
A\in\C^{n\times n}:~H-
{\scriptstyle\frac{1}{{\alpha}^2}}A^*HA\in\overline{\mathbf P}_n\right\}
\end{matrix}\quad\quad\quad H\in\mathbf{H}_n~.
\end{equation}
We now examine the structure of this set.

\begin{Tm}\label{Tm:Stein}
For $H\in\mathbf{H}_n$ the set $\mathbf{S}_H$ 
(${\scriptstyle\frac{1}{\alpha}}\overline{\mathbf S}_H$)
is open (closed) convex closed under multiplication by $c\in\C$,
$1\geq |c|$ and closed under products of its elements.
Furthermore, whenever 
$A\in{\scriptstyle\frac{1}{\alpha}}\mathbf{S}_H$ and
$B\in{\scriptstyle\frac{1}{\beta}}\mathbf{S}_H$, for some
$\alpha$, $\beta>0$, then the product of these matrices satisfies
$AB\in{\scriptstyle\frac{1}{\alpha\beta}}\mathbf{S}_H~$. 
\end{Tm}

{\bf Proof :}
Although classical, for completeness we show this for
${\scriptstyle\frac{1}{\alpha}}\mathbf{S}_H$. Convexity,
closure under multiplication by $c\in\C$, $1\geq |c|$ and the fact
that this set is open, are all trivial.
\vskip 0.2cm

\noindent
Assume that for some $H\in\mathbf{H}_n$ one has that
$A\in{\scriptstyle\frac{1}{\alpha}}\mathbf{S}_H$ and
$B\in{\scriptstyle\frac{1}{\beta}}\mathbf{S}_H$, for some
$\alpha$, $\beta>0$
namely,
\[
\begin{matrix}
H-{\scriptstyle\frac{1}{\alpha^2}}A^*HA=Q_a&&{\rm for~some}~~Q_a\in\mathbf{P}_n
\\~\\
H-{\scriptstyle\frac{1}{\beta^2}}B^*HB=Q_b&&{\rm for~some}~~Q_b\in\mathbf{P}_n~.
\end{matrix}
\]
Multiplying the first equation by ${\scriptstyle\frac{1}{\beta}}B^*$ and
${\scriptstyle\frac{1}{\beta}}B$ from the left and
from the right respectively, and adding the result to the second
equation  yields,
\[
H-{\scriptstyle\frac{1}{\alpha^2\beta^2}}(AB)^*HAB=
{\scriptstyle\frac{1}{\beta^2}}B^*Q_aB+Q_b~,
\]
and as the right hand side is positive definite.
\vskip 0.2cm

\noindent
Establishing this claim for the closed set
${\scriptstyle\frac{1}{\alpha}}\overline{\mathbf S}_H$
is analogous is thus omitted.
\qed
\vskip 0.2cm

\noindent
In the sequel, we focus our attention on the case where in Eq.
\eqref{eq:DefRefinedSteinH} one has that $H\in\mathbf{P}_n~$.

\begin{Cy}\label{Ob:Spectralradius}
Consider the description in Theorem \ref{Tm:Stein} of the set
${\scriptstyle\frac{1}{\alpha}}\mathbf{S}_H$ in Eq.
\eqref{eq:DefRefinedSteinH}. Whenever $H\in\mathbf{P}_n$ this is in addition
a family of matrices whose spectral radius is bounded by $\alpha$.
\end{Cy}

\noindent
Indeed, when $H\in\mathbf{P}_n$, one can multiply the Stein matrix inclusion in Eq.
\eqref{eq:DefRefinedSteinH} by $H^{-\frac{1}{2}}$ from both sides to obtain,
\begin{equation}\label{eq:SteinNorm}
\begin{matrix}
{\scriptstyle\frac{1}{\alpha}}\mathbf{S}_H&
=&\{ A\in\C^{n\times n}:~ \alpha>\| H^{\frac{1}{2}}AH^{-\frac{1}{2}}\|_2\}
\\~\\
{\scriptstyle\frac{1}{\alpha}}\overline{\mathbf{S}}_H&
=&\{ A\in\C^{n\times n}:~\alpha\geq\|
H^{\frac{1}{2}}AH^{-\frac{1}{2}}\|_2\}
\end{matrix}
\quad\quad\quad
\begin{matrix}
H\in\mathbf{P}_n
\\~\\
\alpha>0.
\end{matrix}
\end{equation}
Thus, in particular, the spectral radius of $A$ is bounded by $\alpha$.
\vskip 0.2cm

\noindent
We conclude this section by pointing out that a complete characterization
of the set $\mathbf{S}_H$ in Eq. \eqref{eq:DefSteinH}, for an arbitrary
$H\in\mathbf{H}_n$,
appeared in \cite[Theorem 3.5]{Ando2004}. This remarkable result is
quite involved. Now, on the expense of restricting the case to $H=I_n$,
in Proposition \ref{Pn:MCIC} below, we obtain, through matrix-convexity,
a much simpler characterization. Subsequently, this advantage is
exploited to describe Discrete-time Bounded real rational functions.

\section{Maximal multiplicative matrix-convex sets of matrices}
\label{Sec:MCsets}
\setcounter{equation}{0}

We next resort to the notion of a {\em matrix-convex} set, see e.g.
\cite{EffrWink1997} and more recently, \cite{Evert2018},
\cite{EverHeltKlepMcCull2018}, \cite{HelKlepMcCull2013},
\cite{PassShalSol2018}.

\begin{Dn}\label{Dn:MatrixConvex}
{\rm
A family $\mathbf{A}$, of square matrices (of various dimensions) is
said to be~ {\em matrix-convex}, if for all natural $k$, $n$,
\begin{equation}\label{eq:Isometry1}
\sum\limits_{j=1}^k\upsilon_j^*\upsilon_j=I_n
\quad\quad
\begin{smallmatrix}
\forall{\upsilon}_j\in\C^{{\eta}_j\times n}
\\~\\
\forall\eta_j~,
\end{smallmatrix}
\end{equation}
one has that having $A_1,~\ldots~,~A_k$ (of dimensions
\mbox{${\eta}_1\times{\eta}_1$} through
\mbox{${\eta}_k\times{\eta}_k$}) within
$\mathbf{A}$, implies that also the $n\times n$ matrix
\[
\sum\limits_{j=1}^k\upsilon_j^*A_j\upsilon_j
\]
belongs to $\mathbf{A}$.
\qed
}
\end{Dn}

\noindent
In the sequel, Skew-Hermitian matrices are denoted by, $i\overline{\mathbf H}~$.
It is common to take $\overline{\mathbf H}$ and $i\overline{\mathbf H}$ as the
matricial extension of $\R$ and $i\R$, respectively.
\vskip 0.2cm

\begin{Rk}\label{Rk:MatrixConvex}
{\rm
In \cite{Lewk2020a} it was shown that there are not-too-many, non-trivial
matrix-convex sets, among them:
$\quad\overline{\mathbf H}~,\quad\quad\quad
i\overline{\mathbf H}~,\quad\quad\quad
\overline{\mathbf P},\quad\quad\quad
{\mathbf P}.$
\qed
}
\end{Rk}
\vskip 0.2cm

\noindent
Note that matrix-convexity is rather stringent. Specifically, by definition,
matrix-convexity implies both classical convexity and having the set
invariant under all unitary similarities. The following Example
\ref{Ex:MatrixConvex} illustrates the fact that the converse falls
short from being true.
\vskip 0.2cm

\begin{Ex}\label{Ex:MatrixConvex}
{\rm 
{\bf ~a.} Following Remark \ref{Rk:MatrixConvex}, for arbitrary ${\alpha}>0$,
the~ {\em subset}~ of all matrices in $\overline{\mathbf H}$ (or within
$i\overline{\mathbf H}$ or within $\overline{\mathbf P}$, or within
$\mathbf{P}$) whose spectral radius is less or equal to $\alpha>0$,
is matrix convex.
\vskip 0.2cm

\noindent
{\bf b.} The family of matrices whose Frobenius (a.k.a. Euclidean or
Hilbert-Schmidt) norm, see e.g. \cite[p. 291]{HornJohnson1} is
bounded say by some $\alpha>0$, is both convex and unitarily
invariant, but it is not matrix-convex.
\vskip 0.2cm

\noindent
Consider the set of matrices $\{A~:~5\geq\| A\|_{\rm Frobenius}\}$.
Now from
\mbox{$A=\left(\begin{smallmatrix}4&&0\\~\\0&&3\end{smallmatrix}\right)$}
which belongs to this set $(\| A\|_{\rm Frobenius}=5)$, construct the
matrix
\[
\hat{A}=\underbrace{\left(\begin{smallmatrix}
1&0&0&0\\ 0&0&1&0
\end{smallmatrix}\right)}_{\Upsilon^*}
\left(\begin{matrix}A&0\\0&A\end{matrix}\right)
\underbrace{\left(\begin{smallmatrix}
1&0\\ 0&0\\ 0&1\\ 0&0
\end{smallmatrix}\right)}_{\Upsilon}=4I_2~.
\]
Now, since $\|\hat{A}\|_{\rm Frobenius}=4\sqrt{2}\approx 5.7~$, this set
is not matrix-convex.
}
\qed
\end{Ex}
\vskip 0.2cm

\noindent
We find it convenient to use, for a prescribed $\alpha>0$,
this notation,
\[
{\scriptstyle\frac{1}{\alpha}}\mathbf{S}_I=\bigcup\limits_{n=1}^{\infty}
{\scriptstyle\frac{1}{\alpha}}\mathbf{S}_{I_n}\quad\quad
{\rm and}\quad\quad
{\scriptstyle\frac{1}{\alpha}}\overline{\mathbf S}_I=
\bigcup\limits_{n=1}^{\infty}{\scriptstyle\frac{1}{\alpha}}\overline{
\mathbf S}_{I_n}~.
\]
We next present the key player in this work.

\begin{Pn}\label{Pn:MCIC}
A closed (open), matrix-convex family of matrices whose spectral radius
is less or equal to some $\alpha>0$, containing (on its boundary)
all matrices of the form $\left(\begin{smallmatrix}{\alpha}e^{i\theta}
&&0\\~\\ 0&&0\end{smallmatrix}\right)$ with 
${\scriptstyle\theta}\in[0,~2\pi)$, is the set
${\scriptstyle\frac{1}{\alpha}}\overline{\mathbf S}_I$
(\mbox{${\scriptstyle\frac{1}{\alpha}}\mathbf{S}_I$}).
\vskip 0.2cm

\noindent
Furthermore, the converse is true as well.
\vskip 0.2cm

\noindent
If in addition this is a maximal family of matrices which (whenever
dimensions are suitable) is closed
under products among its elements, then it is equivalent to having
\[
\alpha=1.
\]
\end{Pn}
\vskip 0.2cm

\noindent
{\bf Proof:}\quad
First as matrix-convex sets are in particular closed under unitary
similarity, from the family in the claim, we actually obtain all
$2\times 2$ rank one normal matrices with spectral radius
$\alpha$. Now, taking matrix-convex combinations, yields all
(not necessarily normal) $2\times 2$ matrices with spectral
radius $\alpha$. Taking
convex combination with zero, results in all $2\times 2$
matrices with spectral radius of at most $\alpha$. 
\vskip 0.2cm

\noindent
Let now $U\in\C^{n\times 2}$, where $n\geq 3$, be an arbitrary
isometry, i.e. $U^*U=I_2$. Multiplying a matrix 
\mbox{$A\in\C^{2\times 2}$}
from the above set, by $UAU^*$ yields $n\times n$ matrices 
of degree of at most two, with spectral radius of at most
$\alpha$. Taking further matrix-convex combination
recovers all possible degrees up to (including) $n$.
\vskip 0.2cm

\noindent
To show that this is indeed $\overline{\mathbf S}_I$ (or
$\mathbf{S}_I$) recall that from Eq. \eqref{eq:SteinNorm}
\begin{equation}\label{eq:SteinDisk}
\begin{matrix}
{\scriptstyle\frac{1}{\alpha}}\mathbf{S}_{I_n}&=
\{ A\in\C^{n\times n}:~\alpha>\|A\|_2~\}
&&{\rm and}&&
{\scriptstyle\frac{1}{\alpha}}\overline{\mathbf{S}}_{I_n}&=
\{ A\in\C^{n\times n}:~\alpha\geq\|A\|_2~\}.
\end{matrix}
\end{equation}
Next, recall that for every induced norm, a set of the form
$\{ A\in\C^{n\times n}:~\alpha>\|A\|~\}$, is convex and
the spectral radius of all matrices in it, is bounded
by $\alpha$, see e.g.  \cite[Section 5.6]{HornJohnson1}.
\vskip 0.2cm

\noindent
To guarantee matrix-convexity, the norm must be in addition 
unitarily-invariant, which implies $~\|~\|_2$.
\vskip 0.2cm

\noindent
Next, we show that the closed set
${\scriptstyle\frac{1}{\alpha}}\overline{\mathbf{S}}_{I_n}$ 
(the case of the open set ${\scriptstyle\frac{1}{\alpha}}
\mathbf{S}_{I_n}$ is similar and thus omitted) is matrix-convex.
For a natural parameter $k$ let $\Upsilon\in\C^{kn\times n}$ be an
isometry, i.e. ${\Upsilon}^*\Upsilon=I_n$, then
\[
\begin{matrix}
\left\|{\scriptstyle\Upsilon}^*\left(\begin{smallmatrix}A_1&~&~\\
~&\ddots&~\\~&~&A_k\end{smallmatrix}\right){\scriptstyle\Upsilon}
\right\|_2&\leq&\left\|{\scriptstyle\Upsilon}^*\right\|_2\left\|
\left(\begin{smallmatrix}A_1&~&~\\~&\ddots&~\\~&~&A_k
\end{smallmatrix}\right)\right\|_2\left\|{\scriptstyle\Upsilon}
\right\|_2&&{\rm sub-multiplicative~~norm}\\~\\~&=&\left\|
\left(\begin{smallmatrix}A_1&~&~\\~&\ddots&~\\~&~&A_k
\end{smallmatrix}\right)\right\|_2&&{\scriptstyle\Upsilon}~~{\rm
is~an~isometry}\\~\\~&=&\max\left(\|A_1\|_2~,~\ldots~,~\|A_k\|_2
\right)&&{\rm induced~~norm}\\~\\~&\leq&\alpha&&{\rm assumption},
\end{matrix}
\]
so this part of the claim is established.
\vskip 0.2cm

\noindent
For maximality of the spectral norm under product of elements,
let $B\not\in\overline{\mathbf{S}}_{I_n}$ be arbitrary. One can
always find, within $\mathbf{S}_{I_n}$, a matrix $A$ so that the
spectral radius of the product $AB$, is larger than one (and thus
the spectral radius of $(AB)^l$ is diverging with $l$ natural).
Indeed, let the Singular Value Decomposition, see e.g.
\cite[Theorem 7.35]{HornJohnson1}, of a matrix $B$ be
\[
B=\sum\limits_{j=1}^n{\sigma}_ju_jv_j^*\quad\quad
\begin{smallmatrix}
(1+\epsilon)={\sigma}_1\geq\sigma_2\geq~\cdots~\geq{\sigma}_n\geq 0,
~~~\epsilon>0\\~\\
u_j\in\C^n~~~u_j^*u_k={\delta}_{j,k}~~n\geq j\geq k\geq 1\\~\\
v_j\in\C^n~~~v_j^*v_k={\delta}_{j,k}~~n\geq j\geq k\geq 1,
\end{smallmatrix}
\]
where ${\delta}_{j,k}$ is the Kronecker delta. To avoid triviality, assume
that the spectral radius of $B$ is less than one (Schur stable). This
implies that \mbox{$\frac{1}{1+\epsilon}>|u_1v_1^*|$} (when ${\sigma}_2=0$,
this is in fact sufficient).
\vskip 0.2cm

\noindent
Take now $A={\scriptstyle\frac{1}{1+2\epsilon}}B^*$. By construction
$\| A\|_2={\scriptstyle\frac{1+\epsilon}{1+2\epsilon}}$, so indeed
$A\in\mathbf{S}_{I_n}$. Next,
\[
AB={\scriptstyle\frac{1}{1+2\epsilon}}B^*B={\scriptstyle\frac{1}
{1+2\epsilon}}\sum\limits_{j=1}^n{\sigma}_j^2v_jv_j^*
=
{\scriptstyle\frac{1}{1+2\epsilon}}\left((1+\epsilon)^2v_1v_1^*+
\sum\limits_{j=2}^n{\sigma}_j^2v_jv_j^*\right).
\]
Thus, in fact $AB\in\overline{\mathbf P}_n$ and
$
\| AB\|_2=
{\scriptstyle\frac{(1+\epsilon)^2}{1+2\epsilon}}
=1+
{\scriptstyle\frac{{\epsilon}^2}{1+2\epsilon}},
$
which is also the spectral radius of $AB$, so this
part of the construction is complete.
\vskip 0.2cm

\noindent
The converse direction is to show that the set
${\scriptstyle\frac{1}{\alpha}}\mathbf{S}_{I_n}$ 
is of this structure. This is easy and thus omitted.
\vskip 0.2cm

\noindent
Finally, to obtain a set which is closed under products of
its elements, one needs to take $1\geq\alpha$,
and maximality requires $1=\alpha$. Thus the proof
is complete.
\qed
\vskip 0.2cm

\noindent
As an application consider the following, see e.g.
\cite{MolcPyat1989}.
\vskip 0.2cm

\noindent
{\bf Stability of difference inclusions}
\vskip 0.2cm

\noindent
Recall that the solution $x(j)$ of an autonomous
difference equation $x(j+1)=Ax(j)$ converges to zero for all $x(0)$, if
and only if the spectral radius of $A$ is less than one. Recall also that
the set of matrices whose spectral radius is less than one (colloquially,
``Schur stable") is not closed under multiplication, e.g. both matrices
$A=\left(\begin{smallmatrix}0&2\\0&0\end{smallmatrix}\right)$ and $B=A^*$,
have a zero spectral radius. However, the spectral radius of the product
$AB$, is four. 
\vskip 0.2cm

\noindent
Let $\mathbf{M}$ be a given set of real $n\times n$ matrices.
A difference inclusion
\begin{equation}\label{eq:DifferenceInc}
x(j+1)\in\mathbf{M}x(j)
\quad\quad\quad\begin{smallmatrix}
x(j)\in\R^n\\~\\
j=0,~1,~2,~\ldots
\end{smallmatrix}
\end{equation}
can be interpreted as having
\[
x(j+1)=A(j)x(j)
\quad\quad\quad j=0,~1,~2,~\ldots~~A(j)~~{\rm is}~~
\left\{\begin{matrix}
{\rm arbitrary},
\\~\\
{\rm within}~~\mathbf{M}.
\end{matrix}
\right.
\]
From Proposition \ref{Pn:MCIC} it follows that:

\begin{Cy}
There exists $\alpha\in(0,~1-\epsilon]$, with $1>>\epsilon>0$, so that the
difference inclusion in Eq.  \eqref{eq:DifferenceInc} satisfies Eq.
\eqref{eq:ExpConvergence},
\begin{equation}\label{eq:ExpConvergence}
\|x(0)\|_2{\alpha}^j\geq\|x(j)\|_2\quad 
\forall j=0,~1,~2,~\ldots
\end{equation}
if and only if for the same $\alpha$, 
\[
\mathbf{M}\subset{\scriptstyle\frac{1}{\alpha}}\overline{\mathbf{S}}_{I_n}~.
\]
\end{Cy}

\noindent
For completeness we recall that if for some
\mbox{$\alpha\in(0,~1-\epsilon]$}, with \mbox{$1>>\epsilon>0$}
the condition is relaxed to \mbox{$\mathbf{M}\subset\{ A~:~
\alpha>\| A\|~\}$}, for some induced matrix norm, see e.g.
\cite[Section 5.6]{HornJohnson1}, then
Eq. \eqref{eq:ExpConvergence} holds when
$~\|~\|_2$ is substituted by the above induced norm
\begin{footnote}{In
principle this can further relaxed in two ways: (i) to having possibly another
norm and (ii) ${\scriptstyle\beta}\geq 1$ so that
\mbox{${\scriptstyle\beta}\|x(0)\|{\alpha}^j\geq\|x(j)\|$}
$\forall j=0,~1,~2,~\ldots$}
\end{footnote}.
\vskip 0.2cm

\noindent
In the next section we use the set $\mathbf{S}_{I_n}~$ to
describe a family of rational functions.

\section{Multiplicative Matrix-convex sets of
Rational Functions}
\label{Sec:MCsetsOfSystems}
\setcounter{equation}{0}

In this section we address Discrete-time Bounded real
\mbox{$m\times m$-valued} rational functions $F(z)$,
denoted by $\mathcal{DB}$, satisfying
\begin{equation}\label{eq:MapToUnitDisk1}
\left(I_m-
\left(F(z)\right)^*F(z)\right)
\in\overline{\mathbf P}_m
\quad\quad\quad\forall z\in\C~{\rm s.t.}~|z|>1.
\end{equation}
See e.g. 
\cite{LiuXio2018},
\cite{Najs2013}, \cite{PremJury1994} and \cite{Vaidya1985}.
\vskip 0.2cm

\noindent
Note that Eq. \eqref{eq:MapToUnitDisk1} can be equivalently
written as,
\begin{equation}\label{eq:MapToUnitDisk2}
1\geq\| F(z)\|_2
\quad\quad\quad\forall z\in\C~{\rm s.t.}~|z|>1.
\end{equation}
As already mentioned, $\mathcal{DB}$ functions can be interpreted as the
$Z$-transform of the shift-invariant difference state equation,
\[
x(k+1)=Ax(k)+Bu(k)
\quad\quad
y(k)=Cx(k)+Du(k),
\]
i.e.$~F(z)=C(zI_n-A)^{-1}B+D$, where in the sense of Eq.
\eqref{eq:MapToUnitDisk2}, the input $u$ dominates the output $y$, i.e.
$\|u\|\geq\|y\|$ for all input $u$.
\vskip 0.2cm

\noindent
Eqs. \eqref{eq:MapToUnitDisk1} and \eqref{eq:MapToUnitDisk2}
in particular imply that whenever $F_a$ and $F_b$ are two
\mbox{$m\times m$-valued} rational functions, with this property,
then so is their product $F_aF_b~$.
\vskip 0.2cm

\noindent
Applying the notation of Proposition \ref{Pn:MCIC}, the set of
$\mathcal{DB}$ functions $F(z)$ in Eq. \eqref{eq:MapToUnitDisk2},
can be equivalently
written as $m\times m$-valued rational functions $F(z)$ so
that\begin{footnote}{
Strictly speaking the first line in Eq. \eqref{eq:AlternativeBR}
should be read as saying that:~
``Whenever $z\in\R$ is not a pole of $F((z)$, then
$F(z)\in\R^{m\times m}$."}\end{footnote},
\begin{equation}\label{eq:AlternativeBR}
F(z)\in\left\{\begin{matrix}\R^{m\times m}&z\in\R
\\~\\
\overline{\mathbf{S}}_{I_m}&z\in\C~{\rm s.t.}~|z|>1.
\end{matrix}\right.
\end{equation}
We now illustrate matrix-convexity operations among
matrix-valued rational functions, to be used in the
sequel.

\begin{Ex}\label{Ex:MatrixConvexRational}
{\rm
Let $F_1(z)$, $F_2(z)$ and $F_3(z)$ be rational functions of
dimensions $1\times 1$, $2\times 2$ and $3\times 3$, respectively.
From these functions, by taking matrix-convex operations, one
can construct functions $G_1(z)$, $G_2(z)$ and $G_3(z)$, of
dimensions $1\times 1$, $2\times 2$ and $3\times 3$, respectively.
(To ease reading, the isometric matrices are partitioned
conform ably with $F_1$, $F_2$ and $F_3$):
\[
\begin{matrix}
1=\left(
\begin{smallmatrix}
\frac{6}{7}&0&\frac{2}{7}&0&\frac{3}{7}&0
\end{smallmatrix}\right)
\left(\begin{smallmatrix}
\frac{6}{7}\\0\\ \frac{2}{7}\\0\\\frac{3}{7}\\0
\end{smallmatrix}\right)
&&&
G_1(z)=
\left({\footnotesize\begin{array}{c|cc|ccc}
\frac{6}{7}&0&\frac{2}{7}&0&\frac{3}{7}&0
\end{array}}\right)
\left(\begin{matrix}
F_1(z)&~     &~     \\
~     &F_2(z)&~     \\
~     & ~    &F_3(z)
\end{matrix}\right)
\left({\footnotesize\begin{array}{c}
\frac{6}{7}\\
\hline
0\\ \frac{2}{7}\\
\hline
0\\ \frac{3}{7}\\0
\end{array}}\right),
\end{matrix}
\]
\vskip 0.2cm

\[
\begin{matrix}
I_2=\left(\begin{smallmatrix}0&-\frac{2}{3}&0&\frac{2}{3}&\frac{1}{3}&0\\
\frac{2}{7}&~~\frac{3}{7}&0&0&\frac{6}{7}&0
\end{smallmatrix}\right)
\left(\begin{smallmatrix}
~~0&\frac{2}{7}\\-\frac{2}{3}&\frac{3}{7}\\~~0&0\\~~\frac{2}{3}&0\\~~
\frac{1}{3}&\frac{6}{7}\\~~0&0\end{smallmatrix}\right)
&&
G_2(z)=
\left({\footnotesize\begin{array}{c|cc|ccc}
0&-\frac{2}{3}&0&\frac{2}{3}&\frac{1}{3}&0\\
\frac{2}{7}&~~\frac{3}{7}&0&0&\frac{6}{7}&0
\end{array}}\right)
\left(\begin{matrix}
F_1(z)&~     &~     \\
~     &F_2(z)&~     \\
~     & ~    &F_3(z)
\end{matrix}\right)
\left({\footnotesize\begin{array}{ccc}
~~0&&\frac{2}{7}\\ \hline -\frac{2}{3}&&\frac{3}{7}\\ 
~~0&&0\\ \hline ~~\frac{2}{3}&&0\\
~~\frac{1}{3}&&\frac{6}{7}\\~~0&&0\end{array}}\right),
\end{matrix}
\]
and taking the isometry
$~I_3=\left(\begin{smallmatrix}~~0&0&\frac{3}{5}&0&\frac{4}{5}&0\\
~~\frac{3}{7}&\frac{6}{7}&0&\frac{2}{7}&0&0\\
-\frac{2}{3}&\frac{1}{3}&0&0&0&\frac{2}{3}
\end{smallmatrix}\right)
\left(\begin{smallmatrix}
0&\frac{3}{7}&-\frac{2}{3}\\
0&\frac{6}{7}&~~\frac{1}{3}\\
\frac{3}{5}&0&~~0\\
0&\frac{2}{7}&~~0\\
\frac{4}{5}&0&~~0\\
0&0&~~\frac{2}{3}
\end{smallmatrix}\right)$ yields,
}

\mbox{$G_3(z)=
\left({\footnotesize\begin{array}{c|cc|ccc}
~~0&0&\frac{3}{5}&0&\frac{4}{5}&0\\
~~\frac{3}{7}&\frac{6}{7}&0&\frac{2}{7}&0&0\\
-\frac{2}{3}&\frac{1}{3}&0&0&0&\frac{2}{3}
\end{array}}\right)
\left(\begin{matrix}
F_1(z)&~     &~     \\
~     &F_2(z)&~     \\
~     & ~    &F_3(z)
\end{matrix}\right)
\left({\footnotesize\begin{array}{ccccc}
0&&\frac{3}{7}&&-\frac{2}{3}\\
\hline
0&&\frac{6}{7}&&~~\frac{1}{3}\\
\frac{3}{5}&&0&&~~0\\
\hline
0&&\frac{2}{7}&&~~0\\
\frac{4}{5}&&0&&~~0\\
0&&0&&~~\frac{2}{3}
\end{array}}\right).$}
\qed
\end{Ex}
\vskip 0.2cm

\noindent
Using Eq. \eqref{eq:AlternativeBR} along with Proposition
\ref{Pn:MCIC}, we have the following characterization of 
rational $\mathcal{DB}$ functions.

\begin{Pn}\label{Ob:DTBR}
Let $\mathcal{F}$ be a family of square matrix-valued (of various
dimensions) real rational functions $F(z)$. For all $z$ outside
the closed unit disk, each $F(z)$ is analytic.
\vskip 0.2cm

\noindent
If as a family, $\mathcal{F}$ is matrix-convex\begin{footnote}{In
the sense described in Example \ref{Ex:MatrixConvexRational}.}
\end{footnote} and a maximal set closed under products of its
elements (whenever dimensions are suitable), this is the set
$\mathcal{DB}$ of Discrete-time Bounded real rational functions.
\vskip 0.2cm

\noindent
The converse is true as well.
\end{Pn}

\noindent
As we already mentioned, matrix-convexity is a strong property. This
is next illustrated in the context of $\mathcal{DB}$ functions.

\begin{Ex}\label{Ex:MatrixConvexRationalDB}
{\rm
If $F_1(z)$, $F_2(z)$ and $F_3(z)$ in Example
\ref{Ex:MatrixConvexRational} are $\mathcal{DB}$ functions then
so are $G_1(z)$, $G_2(z)$ and $G_3(z)$.
}
\qed
\end{Ex}
\vskip 0.2cm

\noindent
\begin{Rk}\label{Rk:Rectangular}
{\rm
To simplify the exposition, we consider {\em square} matrix-valued
rational functions. However, some of the results are carried over
to the {\em rectangular} case where,
\[
\overline{\mathbf S}_{I_{m,p}}=\left\{ A\in\C^{p\times m}:~
I_m-A^*A\in\overline{\mathbf P}_m\right\}
\]
and then a $p\times m$-valued rational $\mathcal{DB}$ function
$F(z)$ can be described as\begin{footnote}{
Strictly speaking the first line below
should be read as saying thet:~
``Whenever $z\in\R$ is not a pole of $F((z)$, then
$F(z)\in\R^{p\times m}$."}\end{footnote},
\[
F(z)\in\left\{\begin{matrix}\R^{p\times m}&z\in\R\\~\\
\overline{\mathbf{S}}_{I_{m,p}}&z\in\C~{\rm s.t.}~|z|>1.
\end{matrix}\right.
\]
For more details see e.g. \cite{Najs2013}, and for a subclass
of $\mathcal{DB}$ functions see \cite{AlpJorLew2017}.
\vskip 0.2cm

\noindent
To simplify the exposition, we avoid this generalization.
}
\qed
\end{Rk}
\vskip 0.2cm

\noindent
In the next section we study the structure of families of realization
arrays associated with $\mathcal{DB}$ functions.
\vskip 0.2cm

\section{Sets of Matrix-convex Realization Arrays}
\label{Sec:Realizations}
\setcounter{equation}{0}

\noindent
Recall that whenever $F(z)$ is an $m\times m$-valued rational function
with no pole at infinity, one can associate with it a corresponding
$(n+m)\times(n+m)$ state-space realization array, $R_F$ i.e.
\begin{equation}\label{eq:Realization}
F(z)=C(zI_n-A)^{-1}B+D
\quad\quad\quad
R_F=
\left({\footnotesize
\begin{array}{c|c}
A&B\\
\hline
C&D
\end{array}}\right).
\end{equation}
The realization $R_F$ in Eq. \eqref{eq:Realization} is called minimal,
if $n$ is the McMillan degree of $F(z)$.
\vskip 0.2cm

\noindent
In this section we address {\em families of realization arrays}~ associated
with rational functions. To this end, we adopt the an idea apparently
from \cite[Section 5]{Willems2nd1972} to treat the above
\mbox{$(n+m)\times(n+m)$} $R_F$ as having two faces\begin{footnote}
{Like Janus in the Roman mythology}\end{footnote}: $(i)$ of an
{\em array} and $(ii)$ of a {\em matrix}. (For recent applications of matrix
manipulations of $R_F$ see \cite{AlpayLew2011}, \cite{Lewk2020a},
\cite{Lewk2020b} and \cite{MehrVanDo2020}).
\vskip 0.2cm

\noindent
Before that, a word of caution: For example, $R_1=\left({\footnotesize
\begin{array}{c|c}A&B\\ \hline C&D\end{array}}\right)$ and $R_2=\left(
{\footnotesize\begin{array}{r|r}A&-B\\ \hline-C&D\end{array}}\right)$
are two realization of the same rational function. Furthermore, $R_1$
is minimal (balanced) if and only if $R_2$ is minimal (balanced). However,
$R_3={\scriptstyle\frac{1}{2}}(R_1+R_2)=\left({\footnotesize
\begin{array}{c|c}A&0\\ \hline 0&D\end{array}}\right)$ is a
realization of a {\em zero degree} rational function $F(s)\equiv D$.
\vskip 0.2cm

\noindent
To further study {\em families} of realizations of $\mathcal{DB}$ functions,
we need to introduce a relaxed version of matrix-convexity.

\begin{Dn}\label{Dn:n,mMatrixConvex}
{\rm
For all $k$, let $v_j\in\C^{(n+m)\times(n+m)}$,
$j=1,~\ldots~,~k$ be block-diagonal so that
\begin{equation}\label{eq:n,mIsometry}
\sum\limits_{j=1}^k\underbrace{\left(\begin{smallmatrix}{\upsilon}_{j,n}&0
\\0&{\upsilon}_{j,m}\end{smallmatrix}\right)^*}_{{\upsilon}_j^*}
\underbrace{\left(\begin{smallmatrix}{\upsilon}_{j,n}&0\\0&{\upsilon}_{j,m}
\end{smallmatrix}\right)}_{{\upsilon}_j}
=
\left(\begin{smallmatrix}I_n&&0\\~\\0&&I_m\end{smallmatrix}\right).
\end{equation}
A set $\mathbf{R}$, of $(n+m)\times(n+m)$ matrices, is said to be~
$n,m$-{\em matrix-convex}~ if having
$R_{F_1},~\ldots~,~ R_{F_k}$ in $\mathbf{R}$,
implies that also,
\[
\sum\limits_{j=1}^k
\underbrace{
\left(\begin{smallmatrix}
{\upsilon}_{j,n}&0
\\~\\
0&{\upsilon}_{j,m}
\end{smallmatrix}\right)^*}_{{\upsilon}_j^*}
\underbrace{
\left(\begin{smallmatrix}
A_j&&B_j\\~\\
C_j&&D_j
\end{smallmatrix}\right)}_{R_{F_j}}
\underbrace{
\left(\begin{smallmatrix}
{\upsilon}_{j,n}&0
\\~\\
0&{\upsilon}_{j,m}
\end{smallmatrix}\right)}_{{\upsilon}_j},
\]
belongs to ${\mathbf R}$, for all natural $k$ and all block-diagonal
${\upsilon}_j\in\C^{(n+m)\times(n+m)}$.
\qed
}
\end{Dn}
\vskip 0.2cm

\noindent
In \cite{Lewk2020a} it was pointed out that the notion of
\mbox{$n,m$-{\em matrix-convexity}} is intermediate between (the more strict)
{\em matrix-convexity}, and (weaker) classical convexity.
\vskip 0.2cm

\noindent
For a natural parameter $k$, let  $F_1(z)~,~\ldots~,~F_k(z)$ be a family of
\mbox{$m\times m$-valued} rational functions whose $(n+m)\times(n+m)$
realizations are\begin{footnote}{As poles are all within the unit disk,
these realizations do exist.}\end{footnote},
\begin{equation}\label{eq:RGj}
R_{F_j}=\left({\footnotesize\begin{array}{l|r}\hat{A}_j&\hat{B}_j\\
\hline\hat{C}_j&\hat{D}_j\end{array}}\right)\quad\quad\quad
j=1,~\ldots~,~k.
\end{equation}
Using block-diagonal structured isometries from Eq. \eqref{eq:n,mIsometry}
along with the realizations $R_{F_j}$ in Eq. \eqref{eq:RGj}, let $R_F$
be of the form,
\begin{equation}\label{eq:DefRgBR}
\begin{matrix}
R_F=~\sum\limits_{j=1}^k\underbrace{\left(\begin{smallmatrix}
{\upsilon}_{j,n}&0\\0&{\upsilon}_{j,m}\end{smallmatrix}\right)^*}_{
{\upsilon}_j^*}\underbrace{\left(\begin{smallmatrix}\hat{A}_j&\hat{B}_j\\
\hat{C}_j&\hat{D}_j\end{smallmatrix}\right)}_{R_{F_j}}\underbrace{
\left(\begin{smallmatrix}{\upsilon}_{j,n}&0\\0&{\upsilon}_{j,m}
\end{smallmatrix}\right)}_{{\upsilon}_j}.
\end{matrix}
\end{equation}
Let now $F(z)$ be an \mbox{$m\times m$-valued} rational function whose
realization $R_F$ is given by Eq. \eqref{eq:DefRgBR}. We now address
the following problem: Under what conditions does having the functions
$F_1(z),~\ldots~,~F_k(z)$, in Eq. \eqref{eq:RGj}, Discrete-time-Bounded
real, imply that the resulting $F(z)$ in Eq. \eqref{eq:DefRgBR} is
$\mathcal{DB}$ as well?
\vskip 0.2cm

\noindent
If such a property holds, this suggests that out of a small number of
``extreme points" of balanced realizations of $\mathcal{DB}$ rational
functions, one can construct a whole ``matrix-convex-hull" of
realizations of functions, within the same family. This may enable one
to perform a simultaneous balanced truncation model order reduction of
a whole family of $\mathcal{DB}$ functions, in the spirit of
\cite[Section 5]{CohenLew1997b}.
\vskip 0.2cm

\noindent
As already indicated, even when the ``extreme points" realizations are
balanced, the resulting ``intermediate" realization may be not minimal.
\vskip 0.2cm

\noindent
Recall that the classical version of the Kalman-Yakubovich-Popov
Lemma for Discrete-time Bounded real rational functions, see e.g.
\cite{Najs2013}, \cite{PremJury1994}, \cite{Vaidya1985} and
\cite[Section 4]{XiaoHill1999}, says the following.

\begin{La}\label{La:S}
Let  $F(z)$ be an $m\times m$-valued rational function and let
$R_F$ be a corresponding realization see \eqref{eq:Realization}
\vskip 0.2cm

\noindent
(I)\quad If  there exists a matrix $P\in\mathbf{P}_n$ so that
\begin{equation}\label{eq:OriginalSlemma}
\begin{matrix}
\left(\begin{smallmatrix}P&&0\\~\\0&&I_m\end{smallmatrix}\right)
-
\underbrace{\left(\begin{smallmatrix}A&&B\\~\\C&&D
\end{smallmatrix}\right)^*}_{R_F^*}\left(\begin{smallmatrix}
P&&0\\~\\0&&I_m\end{smallmatrix}\right)\underbrace{
\left(\begin{smallmatrix}A&&B\\~\\C& &D\end{smallmatrix}\right)
}_{R_F}\in\overline{\mathbf{P}}_{n+m}~,
\end{matrix}
\end{equation}
then $F(z)$ is a $\mathcal{DB}$ function.
\vskip 0.2cm

\noindent
If $F(z)$ is in $\mathcal{DB}$ function and its realization in Eq.
\eqref{eq:Realization}, is minimal, i.e. $n$ is the McMillan
degree, then Eq. \eqref{eq:OriginalSlemma} is satisfied.
\vskip 0.2cm

\noindent
(II)\quad Up to change of coordinates, one can substitute 
in Eq. \eqref{eq:OriginalSlemma} $P=I_n$ so that,
\begin{equation}\label{eq:BalancedS}
\left(\begin{smallmatrix}I_n&&0\\~\\0&&I_m\end{smallmatrix}\right)
-
\underbrace{\left(\begin{smallmatrix}A&&B\\~\\C& &D\end{smallmatrix}
\right)^*}_{R_F^*}\left(\begin{smallmatrix}I_n&&0\\~\\0&&I_m
\end{smallmatrix}\right)\underbrace{\left(\begin{smallmatrix}
A&&B\\~\\C& &D\end{smallmatrix}\right)}_{R_F}
\in\overline{\mathbf{P}}_{n+m}~.
\end{equation}
In particular, this is the case when the realization is balanced.
\end{La}
\vskip 0.2cm

\noindent
Combining Proposition \ref{Pn:MCIC} along with part (II) of Lemma
\ref{La:S} we can now answer the question posed in the
beginning of this section.

\begin{Cy}\label{CyConvexRealizationS}
For a natural parameter $k$, let  $F_1(z)~,~\ldots~,~F_k(z)$ be a family of
\mbox{$m\times m$-valued} rational functions. Assume that they all
admit $(n+m)\times(n+m)$ realizations as in Eqs. \eqref{eq:RGj},
\eqref{eq:BalancedS},
i.e.
\begin{equation}\label{eq:SimultaneousSlemma}
\left(\begin{smallmatrix}I_n&&0\\~\\0&&I_m\end{smallmatrix}\right)
-
\underbrace{\left(\begin{smallmatrix}A_j&&B_j\\~\\C_j&&D_j\end{smallmatrix}
\right)^*}_{R_{F_j}^*}\left(\begin{smallmatrix}I_n&&0\\~\\0&&I_m
\end{smallmatrix}\right)\underbrace{\left(\begin{smallmatrix}
A_j&&B_j\\~\\C_j& &D_j\end{smallmatrix}\right)}_{R_{F_j}}
\in\overline{\mathbf{P}}_{n+m}\quad\quad 
\begin{smallmatrix}j=1,~\ldots~,~k\end{smallmatrix}.
\end{equation}
Then, an arbitrary realization $R_F$ defined by Eq. \eqref{eq:DefRgBR},
satisfies Eq. \eqref{eq:BalancedS} and thus the associated $F(z)$, is 
a rational $\mathcal{DB}$ function.
\end{Cy}
\vskip 0.2cm

\noindent
A different formulation (and a different proof) of a similar result
appeared in \cite[Proposition 5.3]{Lewk2020b}.
\vskip 0.2cm

\noindent
The fact that realization arrays of $\mathcal{DB}$ rational functions
are inter-related, is next illustrated.
\vskip 0.2cm

\noindent
\begin{Ex}\label{Ex:MultiplicativeRealizationArrays}
{\rm
We here show how by matrix manipulations of realization arrays of 
Discrete time Bounded rational functions, one can
``generate" a whole family of such functions. 
\vskip 0.2cm

\noindent
{\bf 1.}~ For parameters ${\scriptstyle\theta}\in(0, 1)$ and  $a>1$ 
consider the following scalar $\mathcal{DB}$ rational function of
degree one,
\[
\begin{matrix}
f_1(z)={\scriptstyle\theta}\cdot\frac{a+z}{az+1}&&&R_{f_1}=
\begin{smallmatrix}\frac{1}{a}\end{smallmatrix}\left(
{\footnotesize\begin{array}{c|c}-1&\sqrt{\theta(a^2-1)}\\ \hline
\sqrt{\theta(a^2-1)}&\theta\end{array}}\right).
\end{matrix}
\]
\vskip 0.2cm

\noindent
Treating $R_{f_1}$ as a $2\times 2$ matrix, let us define
\mbox{$R_{f_2}:=
\left(\begin{smallmatrix}1&&~0\\~\\0&&-1
\end{smallmatrix}\right)
R_{f_1}\left(\begin{smallmatrix}-1&&~0\\~\\~0&&~1\end{smallmatrix}\right)$}
so that,
\[
\begin{matrix}
f_2(z)={\scriptstyle\theta}\cdot\frac{a-z}{az-1}&&&R_{f_2}=
\begin{smallmatrix}\frac{1}{a}\end{smallmatrix}\left(
{\footnotesize\begin{array}{c|c}1&
\sqrt{\theta(a^2-1)}\\ \hline
\sqrt{\theta(a^2-1)}&-\theta\end{array}}\right).
\end{matrix}
\]
is another Discrete-time Bounded scalar rational functions of degree one.
\vskip 0.2cm

\noindent
For example, taking now
$R_{f_3}:=\frac{1}{2}\left(R_{f_1}+R_{f_2}\right)$ yields yet
another $\mathcal{DB}$ rational functions of degree one,
\[
\begin{matrix}
f_3(z)={\scriptstyle\frac{\theta}{a^2}}\cdot\frac{a^2-1}{z}&&&R_{f_3}=\begin
{smallmatrix}\sqrt{\theta(1-\frac{1}{a^2})}\end{smallmatrix}\left(
{\footnotesize\begin{array}{c|c}0&1\\ \hline 1&0\end{array}}\right).
\end{matrix}
\]
Note that each of the realization $R_{f_1}$, $R_{f_2}$ and
$R_{f_3}$ is balanced and satisfies Eq. \eqref{eq:SimultaneousSlemma}.
\vskip 0.2cm

\noindent
The following product, $f_4(z):=f_1(z)f_2(z)$, is a $\mathcal{DB}$ rational
function of degree two,
\begin{equation}\label{eq:F5}
\begin{matrix}
f_4(z)={\scriptstyle\theta}^2\cdot\frac{a^2-z^2}{a^2z^2-1}
&&&
R_{f_4}=\begin{smallmatrix}\frac{1}{a}\end{smallmatrix}\left(
{\footnotesize\begin{array}{cc|c}
-1&\frac{\theta}{a}(1-a^2)&
\frac{\theta}{a}\sqrt{\theta(a^2-1)}\\
0&1&\sqrt{\theta(a^2-1)}\\ \hline
-\sqrt{\theta(a^2-1)}&\frac{\theta}{a}\sqrt{\theta(a^2-1)}
&-\frac{{\theta}^2}{a}
\end{array}}\right).
\end{matrix}
\end{equation}
Now, 
\mbox{$\hat{R}_{f_4}:=\left(\begin{smallmatrix}0&-1&~0\\1&~0&~0\\0&~0&~1
\end{smallmatrix}\right)R_{f_4}\left(\begin{smallmatrix}~~0&~1&~0\\
-1&~0&~0\\~~0&~0&~1\end{smallmatrix}\right)$} is another minimal
realization of the same $f_4(z)$. Let now define
$R_{f_5}:=\frac{1}{2}\left(R_{f_4}+\hat{R}_{f_4}\right)$.
\vskip 0.2cm

\noindent
For simplicity of presentation take $\theta=\frac{1}{2}$ and
$a=3$ then, out of $f_4(z)$, whose poles are
$\pm{\scriptstyle\frac{1}{3}}$, one
obtains the following $\mathcal{DB}$ function, of degree two,
\[
\begin{matrix}
f_5(z)=\frac{1}{36}\cdot\frac{\left(\frac{4}{3}\right)^2-z^2}
{z^2+\left(\frac{2}{9}\right)^2}&&&R_{f_5}=\begin{smallmatrix}
\frac{1}{36}\end{smallmatrix}\left({\footnotesize\begin{array}
{cr|r}0&-8&-10\\8&0&14\\ \hline 14&10&1\end{array}}\right),
\end{matrix}
\]
whose poles are $\pm{\scriptstyle\frac{2}{9}}i$.
\vskip 0.2cm

\noindent
{\bf 2.}~ In a way similar to the above, 
we next illustrate how by taking matrix products and matrix-convex
operations, one can ``generate"} virtually all {\rm realization arrays
of a scalar $\mathcal{DB}$ rational of McMillan degree of at most 1.
\vskip 0.2cm

\noindent
Indeed, substituting in Eq. \eqref{eq:BalancedS} $n=1$ and $m=1$, 
reveals that one can take $R_{f_j}$ ($j$ natural) to be a $2\times 2$
orthogonal matrix (with ${\rm det}=1$) of the form
\begin{equation}\label{eq:2x2orthogonal}
f_j(z)=\begin{smallmatrix}\frac{z\cos({\theta}_j)-1}{z-\cos({\theta}_j)}
\end{smallmatrix}\quad\quad R_{f_j}:=
\left(\begin{smallmatrix}\cos({\theta}_j)&&-\sin({\theta}_j)\\~\\
\sin({\theta}_j)&&~~\cos({\theta}_j)\end{smallmatrix}\right)
\quad\quad\quad
\begin{smallmatrix}j=1,~2,~\ldots
\\~\\
{\theta}_j\in[0,~2\pi).\end{smallmatrix}
\end{equation}
Recall also that for all ${\theta}_1$ and ${\theta}_2$,
\[
R_{f_1}R_{f_2}=
\left(\begin{smallmatrix}\cos({\theta}_1+{\theta}_2)&&
-\sin({\theta}_1+{\theta}_2)\\~\\
\sin({\theta}_1+{\theta}_2)&&~~
\cos({\theta}_1+{\theta}_2)\end{smallmatrix}\right),
\]
which can be viewed as both: another orthogonal matrix of the form
of Eq. \eqref{eq:2x2orthogonal} and a balanced realization another
$\mathcal{DB}$ function satisfying Eq. \eqref{eq:BalancedS} with
zero right hand side, colloquially para-unitary, see e.g.
\cite{AlpJorLew2017}.
\vskip 0.2cm

\noindent
For simplicity, take now an arbitrary \mbox{${\theta}_1={\theta}_2$}
provided that $\frac{{\theta}_1}{\pi}$ is irrational (e.g. ${\theta}_1$
rational). This means that by taking infinite powers of the form
$(R_{f_1})^k$, $k=1,~2,~\ldots~$ one obtains a dense subset of all
$2\times 2$ orthogonal matrices of the from of Eq.
\eqref{eq:2x2orthogonal}.
\vskip 0.2cm

\noindent
Next take the matrix \mbox{$\tilde{R}=\left(\begin{smallmatrix}-1&~~0
\\~~0&~~1\end{smallmatrix}\right)$}. $\tilde{R}$ can also be viewed
as a (non-minimal) realization of the zero degree
rational function $f(z)\equiv 1$.\quad As matrices, we have,
\[
\underbrace{
\left(\begin{smallmatrix}-1&&0\\~\\~~0&&1\end{smallmatrix}\right)
}_{\tilde{R}}
\underbrace{
\left(\begin{smallmatrix}\cos({\theta}_j)&&-\sin({\theta}_j)\\~\\
\sin({\theta}_j)&&~~\cos({\theta}_j)\end{smallmatrix}\right)
}_{R_{f_j}}
=
\underbrace{
\left(\begin{smallmatrix}-\cos({\theta}_j)&&\sin({\theta}_j)\\~\\
~~\sin({\theta}_j)&&\cos({\theta}_j)\end{smallmatrix}\right)
}_{R_{\hat{f}_j}}.
\]
For each $j=1,~2,~\ldots$ the right hand side, $R_{\hat{f}_j}$ is a
symmetric orthogonal matrix (with ${\rm det}=-1$) and a realization
of,
\[
\hat{f}_j(z)
=\begin{smallmatrix}\frac{z\cos({\theta}_j)+1}{z+\cos({\theta}_j)}
\end{smallmatrix}\quad\quad\quad
\begin{smallmatrix}j=1,~2,~\ldots\end{smallmatrix}
\]
Finally, note that $f_j(z)$ and $\hat{f}_j(z)$ are weak contractions.
Taking matrix-convex combinations of $R_{f_j}$ and of $R_{\hat{f}_j}$,
yields strict contractions. Thus,
the sought construction of all arrays is complete.
\vskip 0.2cm

\noindent
{\bf 3.}~ Recall that in principle the idea of item 2 can be carried
over to higher dimensions, by taking products of planar rotations
(a.k.a. Givens rotations, see e.g.
\cite[Example 2.2.3]{HornJohnson1}) of the form,
\vskip 0.2cm

$
\left(\begin{smallmatrix}\cos({\theta}_1)&&-\sin({\theta}_1)&&0\\~\\
\sin({\theta}_1)&&~~\cos({\theta}_1)&&0\\~\\
0&&0&&1\end{smallmatrix}\right)
\left(\begin{smallmatrix}\cos({\theta}_2)&&0&&-\sin({\theta}_2)\\~\\
0&&1&&~0\\~\\
\sin({\theta}_2)&&0&&~~\cos({\theta}_2)\end{smallmatrix}\right)
\left(\begin{smallmatrix}1&&0&&~0\\~\\0&&\cos({\theta}_3)&&
-\sin({\theta}_3)\\~\\
0&&\sin({\theta}_3)&&~~\cos({\theta}_3)\end{smallmatrix}\right)
\quad
\quad
\quad
\begin{smallmatrix}
{\theta}_1,
{\theta}_2,
{\theta}_3\in[0,~2\pi).
\end{smallmatrix}
$
\vskip 0.2cm

\noindent
{\bf 4.}~ We conclude by pointing out that one can be more adventurous
in manipulating realization arrays. Take for instance $R_{f_4}$ from
Eq. \eqref{eq:F5}. By construction, as a matrix it satisfies
\[
\left(I_3-{R_{f_4}}^*R_{f_4}\right)\in\overline{\mathbf P}_3~.
\]
Now, another partitioning of $R_{f_4}$ corresponds to a
$2\times 2$-valued
$\mathcal{DB}$ function which is of McMillan degree one,
\[
\begin{matrix}
F_6(z)=\begin{smallmatrix}\frac{1}{a}\end{smallmatrix}\left(
\begin{smallmatrix}1&&\sqrt{\theta(a^2-1)}\\~\\
\theta\sqrt{\theta(a^2-1)}\cdot\frac{z+a}{az+1}
&&-\frac{{\theta}^2}{a}\cdot\frac{z+a}{az+1}\end{smallmatrix}\right)
&&
R_{F_6}=\begin{smallmatrix}\frac{1}{a}\end{smallmatrix}\left(
{\footnotesize\begin{array}{c|cc}-1&\frac{\theta}{a}(1-a^2)&
\frac{\theta}{a}\sqrt{\theta(a^2-1)}\\
\hline
0&1&\sqrt{\theta(a^2-1)}\\
-\sqrt{\theta(a^2-1)}&\frac{\theta}{a}\sqrt{\theta(a^2-1)}&
-\frac{{\theta}^2}{a}
\end{array}}\right).
\end{matrix}
\]
Finally, clearly \mbox{$
\left(\begin{smallmatrix}0&&1\end{smallmatrix}\right)
F_6(z)
\left(\begin{smallmatrix}0\\~\\1\end{smallmatrix}\right)$}
is a scalar $\mathcal{DB}$ function of McMillan degree one.
In fact it equals to \mbox{$~-\frac{\theta}{a}\cdot f_1(z)
=-\frac{{\theta}^2}{a^2}\cdot\frac{z+a}{z+\frac{1}{a}}~$}.
}
\qed
\end{Ex}
\vskip 0.2cm

\begin{center}
ACKNOWLEDGEMENT
\end{center}
\vskip 0.2cm

\noindent
The constructive and exceptionally thorough review, is highly appreciated.

\end{document}